\documentclass[12pt,a4paper]{article}
\usepackage{amsmath}
\usepackage{amssymb}
\usepackage{enumerate}
\usepackage{t1enc}
\usepackage[latin1]{inputenc}
\usepackage[english]{babel}
\usepackage{enumerate}

\usepackage{mathtools}

\usepackage{amsfonts}
\usepackage{latexsym}
\usepackage{bm}
\newtheorem{theorem}{Theorem}[section]

\newtheorem{lemma}[theorem]{Lemma}
\newtheorem{cor}[theorem]{Corollary}

\begin{document}
\title{Isolation of cycles}

\author{Peter Borg\\[5mm]
{\normalsize Department of Mathematics} \\
{\normalsize Faculty of Science} \\
{\normalsize University of Malta}\\
{\normalsize Malta}\\
{\normalsize \texttt{peter.borg@um.edu.mt}}
}

\date{}
\maketitle

\begin{abstract}
For any graph $G$, let $\iota_{\rm c}(G)$ denote the size of a smallest set $D$ of vertices of $G$ such that the graph obtained from $G$ by deleting the closed neighbourhood of $D$ contains no cycle. We prove that if $G$ is a connected $n$-vertex graph that is not a triangle, then $\iota_{\rm c}(G) \leq n/4$. We also show that the bound is sharp. Consequently, we solve a problem of Caro and Hansberg.
\end{abstract}

\section{Introduction}
Unless stated otherwise, we use small letters such as $x$ to denote non-negative integers or elements of sets, and capital letters such as $X$ to denote sets or graphs. The set of positive integers is denoted by $\mathbb{N}$. For $n \geq 1$, $[n]$ denotes the set $\{1, \dots, n\}$ (that is, $[n]= \{i \in \mathbb{N} \colon i \leq n\}$). We take $[0]$ to be the empty set $\emptyset$. Arbitrary sets are assumed to be finite. For a set $X$, ${X \choose 2}$ denotes the set of $2$-element subsets of $X$ (that is, ${X \choose 2} = \{ \{x,y \} \colon x,y \in X, x \neq y \}$).

If $Y$ is a subset of ${X \choose 2}$ and $G$ is the pair $(X,Y)$, then $G$ is called a \emph{graph}, $X$ is called the \emph{vertex set of $G$} and is denoted by $V(G)$, and $Y$ is called the \emph{edge set of $G$} and is denoted by $E(G)$. A \emph{vertex of $G$} is an element of $V(G)$, and an \emph{edge of $G$} is an element of $E(G)$. We call $G$ an \emph{$n$-vertex graph} if $|V(G)| = n$. We may represent an edge $\{v,w\}$ by $vw$. If $vw \in E(G)$, then we say that $w$ is a \emph{neighbour of $v$ in $G$} (and vice-versa). For $v \in V(G)$, $N_{G}(v)$ denotes the set of neighbours of $v$ in $G$, $N_{G}[v]$ denotes $N_{G}(v) \cup \{ v \}$, 
and $d_{G}(v)$ denotes $|N_{G} (v)|$ 
and is called the \emph{degree of $v$ in $G$}. 
For $S \subseteq V(G)$, $N_G[S]$ denotes $\bigcup_{v \in S} N_G[v]$ (the \emph{closed neighbourhood of $S$}), $G[S]$ denotes $(S,E(G) \cap {S \choose 2})$ (the \emph{subgraph of $G$ induced by $S$}), and $G-S$ denotes $G[V(G) \backslash S]$ (the graph obtained by \emph{deleting $S$} from $G$). Where no confusion arises, the subscript $G$ may be omitted from the notation above that uses it; for example, $N_G(v)$ may be abbreviated to $N(v)$. 

If $G$ and $H$ are graphs, $f : V(H) \rightarrow V(G)$ is a bijection, and $E(G) = \{f(v)f(w)\,$: $vw \in E(H)\}$, then we say that $G$ is a \emph{copy of $H$}, and we write $G \simeq H$. Thus, a copy of $H$ is a graph obtained by relabeling the vertices of $H$.
 
For $n \geq 1$, the graphs $([n], {[n] \choose 2})$ and $([n], \{\{i,i+1\} \colon i \in [n-1]\})$ are denoted by $K_n$ and $P_n$, respectively. For $n \geq 3$, $C_n$ denotes the graph $([n], \{\{1,2\}, \{2,3\}, \dots, \{n-1,n\}, \{n,1\}\})$ ($= ([n], E(P_n) \cup \{n,1\})$). A copy of $K_n$ is called a \emph{complete graph}. 
A copy of $P_n$ is called an \emph{$n$-path} or simply a \emph{path}. A copy of $C_n$ is called an \emph{$n$-cycle} or simply a \emph{cycle}. We call a $3$-cycle a \emph{triangle}. Note that $K_3$ is the triangle $C_3$.

If $G$ and $H$ are graphs such that $V(H) \subseteq V(G)$ and $E(H) \subseteq E(G)$, then $H$ is called a \emph{subgraph of $G$}, and we say that \emph{$G$ contains $H$}. 

If $\mathcal{F}$ is a set of graphs and $F$ is a copy of a graph in $\mathcal{F}$, then we call $F$ an \emph{$\mathcal{F}$-graph}. If $G$ is a graph and $D \subseteq V(G)$ such that $G-N[D]$ contains no $\mathcal{F}$-graph, then $D$ is called an \emph{$\mathcal{F}$-isolating set of $G$}. Let $\iota(G, \mathcal{F})$ denote the size of a smallest $\mathcal{F}$-isolating set of $G$. We abbreviate $\iota(G, \{F\})$ to $\iota(G, F)$. 
The study of isolating sets was introduced by Caro and Hansberg~\cite{CaHa17}. It is an appealing and natural generalization of the classical domination problem \cite{C, CH, HHS, HHS2, HL, HL2}. Indeed, $D$ is a $\{K_1\}$-isolating set of $G$ if and only if $D$ is a \emph{dominating set of $G$} (that is, $N[D] = V(G)$), so the $\{K_1\}$-isolation number is the \emph{domination number} (the size of a smallest dominating set). Let $\mathcal{C}$ denote $\{C_k \colon k \geq 3\}$. In this paper, we obtain a sharp upper bound for $\iota(G, \mathcal{C})$, and consequently we solve a problem of Caro and Hansberg \cite{CaHa17}. 

We call a subset $D$ of $V(G)$ a \emph{cycle isolating set of $G$} if $G-N[D]$ contains no cycle (that is, $G-N[D]$ is a \emph{forest}). We denote the size of a smallest cycle isolating set of $G$ by $\iota_{\rm c}(G)$. Thus, $\iota_{\rm c}(G) = \iota(G, \mathcal{C})$.

If $G_1, \dots, G_t$ are graphs such that $V(G_i) \cap V(G_j) = \emptyset$ for every $i,j \in [t]$ with $i \neq j$, then $G_1, \dots, G_t$ are \emph{vertex-disjoint}. A graph $G$ is \emph{connected} if, for every $v, w \in V(G)$, $G$ contains a path $P$ with $v, w \in V(P)$. A connected subgraph $H$ of $G$ is a \emph{component of $G$} if, for each connected subgraph $K$ of $G$ with $K \neq H$, $H$ is not a subgraph of $K$. Clearly, any two distinct components of $G$ are vertex-disjoint. 

For $n, k \in \mathbb{N}$, let $a_{n,k} = \left \lfloor \frac{n}{k+1} \right \rfloor$ and $b_{n,k} = n - ka_{n,k}$. Thus, $a_{n,k} \leq b_{n,k} \leq a_{n,k} + k$. If $F$ is a $k$-vertex graph and $n \leq k$, then let $B_{n,F} = P_n$. If $F$ is a $k$-vertex graph and $n \geq k+1$, then let $F_1, \dots, F_{a_{n,k}}$ be copies of $F$ such that $P_{b_{n,k}}, F_1, \dots, F_{a_{n,k}}$ are vertex-disjoint, let $P_{b_{n,k}}^* = ([b_{n,k}], E(P_{a_{n,k}}) \cup \{a_{n,k}i \colon i \in [b_{n,k}] \backslash [a_{n,k}]\})$, and let $B_{n,F}$ be the connected $n$-vertex graph given by 
\[B_{n,F} = \left( V(P_{b_{n,k}}^* ) \cup \bigcup_{i=1}^{a_{n,k}} V(F_i), \, E(P_{b_{n,k}}^*) \cup \{iv \colon i \in [a_{n,k}], v \in V(F_i)\} \cup \bigcup_{i=1}^{a_{n,k}} E(F_i) \right).\] 
Thus, $B_{n,F}$ is the graph obtained by taking $P_{b_{n,k}}^*, F_1, \dots, F_{a_{n,k}}$ and joining $i$ (a vertex of $P_{b_{n,k}}^*$) to each vertex of $F_i$ for each $i \in [a_{n,k}]$.

For any $n \in \mathbb{N}$, any family $\mathcal{F}$ of graphs, and any $F \in \mathcal{F}$, let 
\[\iota(n,\mathcal{F},F) = \max \left\{\iota(G,\mathcal{F}) \colon G \mbox{ is a connected graph}, V(G) = [n], G \not\simeq F \right\}.\]
We abbreviate $\iota(n,\{F\},F)$ to $\iota(n,F)$. Let $\iota_{\rm c}(n) = \iota(n,\mathcal{C},K_3)$. In Section~\ref{Proofsection}, we prove the following result.

\begin{theorem} \label{result}
If $G$ is a connected $n$-vertex graph that is not a triangle, then
\[\iota_{\rm c}(G) \leq \frac{n}{4}.\] 
Consequently, for any $n \geq 1$, 
\[\iota_{\rm c}(n) = \iota_{\rm c}(B_{n,K_3}) = \left \lfloor \frac{n}{4} \right \rfloor .\]
\end{theorem}
The equality $\iota_{\rm c}(B_{n,K_3}) = \left \lfloor \frac{n}{4} \right \rfloor$ is generalized in the following result.

\begin{lemma}\label{thmlemma} Let $n, k \in \mathbb{N}$ and let $F$ be a $k$-vertex graph.\medskip 
\\
(i) 
If $n \neq k$ or $F \not\simeq P_k$, then
\[\iota(B_{n,F}, F) = \left \lfloor \frac{n}{k+1} \right \rfloor.\]
(ii) If $\mathcal{F}$ is a family of graphs, $F \in \mathcal{F}$, and $n \neq k = |V(F)| = \min\{|V(H)| \colon H \in \mathcal{F}\}$, then
\[\left \lfloor \frac{n}{k+1} \right \rfloor = \iota(B_{n,F}, F) \leq \iota(n,F) \leq \iota(n,\mathcal{F},F).\]
\end{lemma}
\textbf{Proof.} Let $B = B_{n,F}$. If either $n < k$ or $n = k$ and $F \not\simeq P_k$, then $\iota(B, F) = 0$. Suppose $n \geq k+1$. Then, $\iota(B, F) \leq a_{n,k}$ as $[a_{n,k}]$ is a dominating set of $B$. Let $D$ be an $\{F\}$-isolating set of $B$ of size $\iota(B, F)$. For each $i \in [a_{n,k}]$, $D \cap (V(F_i) \cup \{i\}) \neq \emptyset$ as $B - N_{B}[D]$ does not contain the copy $F_i$ of $F$. Thus, $|D| \geq a_{n,k}$. Hence, (i) is proved.

Let $\mathcal{F}$ and $n$ be as in (ii). Since $n \neq |V(F)|$, $B \not\simeq F$. Since $|V(B)| = n$, we can choose a copy $B'$ of $B$ with $V(B') = [n]$. Since $B$ is connected, $B'$ is connected. Thus, $\iota(B',F) \leq \iota(n,F)$. Since $F \in \mathcal{F}$, the $\mathcal{F}$-isolating sets of a graph $G$ are $\{F\}$-isolating sets of $G$, so $\iota(G,F) \leq \iota(G,\mathcal{F})$. Thus, $\iota(n,F) \leq \iota(n,\mathcal{F},F)$. Now $\iota(n,B') = \iota(n,B) = \left \lfloor \frac{n}{k+1} \right \rfloor$ by (i). Hence, (ii) is proved.~\hfill{$\Box$} 
\\

By the results above, $\left \lfloor \frac{n}{4} \right \rfloor$ is a sharp upper bound on $\iota(G,K_3)$ for connected $n$-vertex graphs $G \not\simeq K_3$.

\begin{theorem} For any $n \geq 1$, 
\[\iota(n,K_3) = \iota(B_{n,K_3}, K_3) = \left \lfloor \frac{n}{4} \right \rfloor .\]
\end{theorem}
\textbf{Proof.} Let $G$ be a connected $n$-vertex graph that is not a copy of $K_3$. By Lemma~\ref{thmlemma} (ii), $\left \lfloor \frac{n}{4} \right \rfloor = \iota(B_{n,K_3}, K_3)\leq \iota(n,K_3) \leq \iota(n,\mathcal{C},K_3)$ as $K_3 = C_3 \in \mathcal{C}$. By Theorem~\ref{result}, $\iota(n,\mathcal{C},K_3) = \left \lfloor \frac{n}{4} \right \rfloor$. The result follows.~\hfill{$\Box$}\\

In \cite{CaHa17}, Caro and Hansberg showed that $\frac{1}{4} \leq \limsup_{n \rightarrow \infty} \frac{\iota_{\rm c}(n)}{n} \leq \frac{1}{3}$. In Problem~7.3 of the same paper, they asked for the value of $\limsup_{n \rightarrow \infty} \frac{\iota_{\rm c}(n)}{n}$. The answer is immediately given by Theorem~\ref{result}.

\begin{cor} $\limsup_{n \rightarrow \infty} \frac{\iota_{\rm c}(n)}{n} = \frac{1}{4}$.
\end{cor}
\textbf{Proof.} By Theorem~\ref{result}, for any $n \in \mathbb{N}$, we have $\frac{1}{4} - \frac{3}{4n} = \frac{1}{n}\left( \frac{n-3}{4} \right) \leq \frac{\iota_{\rm c}(n)}{n} \leq \frac{1}{4}$, and, if $n$ is a multiple of $4$, then $\frac{\iota_{\rm c}(n)}{n} = \frac{1}{4}$. Thus, $\lim_{n \rightarrow \infty} \sup \left\{\frac{\iota_{\rm c}(k)}{k} \colon k \geq n \right\}  = \lim_{n \rightarrow \infty} \frac{1}{4} = \frac{1}{4}$.~\hfill{$\Box$}

\section{Proof of Theorem~\ref{result}} \label{Proofsection}

In this section, we prove Theorem~\ref{result}. We start with two lemmas that will be used repeatedly. 

\begin{lemma} \label{lemma}
If $G$ is a graph, $\mathcal{F}$ is a set of graphs, $X \subseteq V(G)$, and $Y \subseteq N[X]$, then \[\iota(G, \mathcal{F}) \leq |X| + \iota(G-Y, \mathcal{F}).\] 
\end{lemma}
\textbf{Proof.} Let $D$ be an $\mathcal{F}$-isolating set of $G-Y$ of size $\iota(G-Y, \mathcal{F})$. Clearly, $\emptyset \neq V(F) \cap Y \subseteq V(F) \cap N[X]$ for each $\mathcal{F}$-graph $F$ that is a subgraph of $G$ and not a subgraph of $G-Y$. Thus, $D \cup X$ is an $\mathcal{F}$-isolating set of $G$. The result follows.~\hfill{$\Box$}
\\

For a graph $G$ and a set $\mathcal{F}$ of graphs, let ${\rm C}(G)$ denote the set of components of $G$, and let ${\rm C}(G, \mathcal{F}) = \{H \in {\rm C}(G) \colon H \mbox{ is an $\mathcal{F}$-graph}\}$. We abbreviate ${\rm C}(G, \{K_3\})$ to ${\rm C}'(G)$. Thus, ${\rm C}'(G) = \{H \in {\rm C}(G) \colon H \mbox{ is a triangle}\}$.

\begin{lemma} \label{lemmacomp}
If $G$ is a graph and $\mathcal{F}$ is a set of graphs, then \[\iota(G, \mathcal{F}) = \sum_{H \in {\rm C}(G)} \iota(H, \mathcal{F}).\] 
\end{lemma}
\textbf{Proof.} For each $H \in {\rm C}(G)$, let $D_H$ be a smallest $\mathcal{F}$-isolating set of $H$. Then, $\bigcup_{H \in {\rm C}(G)} D_H$ is an $\mathcal{F}$-isolating set of $G$, so $\iota(G,\mathcal{F}) \leq \sum_{H \in {\rm C}(G)} |D_H| = \sum_{H \in {\rm C}(G)} \iota(H,\mathcal{F})$. Let $D$ be a smallest $\mathcal{F}$-isolating set of $G$. For each $H \in {\rm C}(G)$, $D \cap V(H)$ is an $\mathcal{F}$-isolating set of $H$. We have $\sum_{H \in {\rm C}(G)} \iota(H,\mathcal{F}) \leq \sum_{H \in {\rm C}(G)} |D \cap V(H)| = |D| = \iota(G,\mathcal{F})$. The result follows.~\hfill{$\Box$}
\\
\\
\textbf{Proof of Theorem~\ref{result}.} Let us first assume the bound 
in the first part of the theorem. Then, $\iota_{\rm c}(n) \leq \frac{n}{4}$. Since $\iota_{\rm c}(n)$ is an integer, $\iota_{\rm c}(n) \leq \left \lfloor \frac{n}{4} \right \rfloor$. Together with Lemma~\ref{thmlemma} (i), this gives us $\iota_{\rm c}(n) = \iota_{\rm c}(B_{n,K_3}, K_3) = \left \lfloor \frac{n}{4} \right \rfloor$.

We now prove the first part of the theorem. We use induction on $n$.   Let $G$ be a connected $n$-vertex graph that is not a triangle. If $n \leq 3$, then, since $G$ is not a triangle, $G$ contains no cycle, and hence $\iota_{\rm c}(G) = 0$. Suppose $n \geq 4$. Let $k = \max \{d(v) \colon v \in V(G)\}$. Since $G$ is connected, $k \geq 2$. Let $v \in V(G)$ such that $d(v) = k$. If $k = 2$, then $G$ is a path or a cycle, $\{v\}$ is a cycle isolating set of $G$, and hence $\iota_{\rm c}(G) \leq 1 \leq \frac{n}{4}$. Suppose $d(v) \geq 3$. Then, $|N[v]| \geq 4$. If $V(G) = N[v]$, then $\{v\}$ is a cycle isolating set of $G$, so $\iota(G) \leq 1 \leq \frac{n}{4}$. Suppose $V(G) \neq N[v]$. Let $G' = G-N[v]$ and $n' = |V(G')|$. Then, $n \geq n' + 4$ and $V(G') \neq \emptyset$. Let $\mathcal{H} = {\rm C}(G')$ and $\mathcal{H}' = {\rm C}'(G')$. By the induction hypothesis, 
\begin{equation} \iota_{\rm c}(H) \leq \frac{|V(H)|}{4} \quad \mbox{for each } H \in \mathcal{H} \backslash \mathcal{H}'. \nonumber
\end{equation}
If $\mathcal{H}' = \emptyset$, then, by Lemma~\ref{lemma} (with $X = \{v\}$ and $Y = N[v]$) and Lemma~\ref{lemmacomp},
\begin{equation}
\iota_{\rm c}(G) \leq 1 + \iota_{\rm c}(G') = 1 + \sum_{H \in \mathcal{H}} \iota_{\rm c}(H) \leq 1 + \sum_{H \in \mathcal{H}} \frac{|V(H)|}{4} = \frac{4 + n'}{4} \leq \frac{n}{4}. \nonumber
\end{equation}

Suppose $\mathcal{H}' \neq \emptyset$. For any $H \in \mathcal{H}$ and any $x \in N(v)$ such that $xy \in E(G)$ for some $y \in V(H)$, we say that $H$ is \emph{linked to $x$} and that $x$ is \emph{linked to $H$}. Since $G$ is connected, each member of $\mathcal{H}$ is linked to at least one member of $N(v)$. Let $L = \{x \in N(v) \colon x \mbox{ is linked to some member of } \mathcal{H}'\}$. Since $\mathcal{H}' \neq \emptyset$, $L \neq \emptyset$. Let $x \in L$. Let $\mathcal{H}'_x = \{H \in \mathcal{H}' \colon H \mbox{ is linked to } x\}$ and $\mathcal{H}_x^* = \{H \in \mathcal{H} \backslash \mathcal{H}' \colon H \mbox{ is linked to $x$ only}\}$. Let $U = \bigcup_{H \in \mathcal{H}'} V(H)$. Let $U_x = N(x) \cap U$ and $U_x^+ = \{x\} \cup U_x$. Note that if a component $A$ of $G - U_x^+$ is a triangle, then $V(A) = N[v] \backslash \{x\}$. 

Suppose $|U_x| \geq 3$. If no component of $G - U_x^+$ is a triangle, then, by Lemma~\ref{lemma} (with $X = \{x\}$ and $Y = U_x^+$), Lemma~\ref{lemmacomp}, and the induction hypothesis, we have 
\begin{equation} \iota_{\rm c}(G) \leq 1 + \iota_{\rm c}(G - U_x^+) = 1 + \sum_{H \in {\rm C}(G - U_x^+)} \iota_{\rm c}(H) \leq \frac{|U_x^+|}{4} + \sum_{H \in {\rm C}(G - U_x^+)} \frac{|V(H)|}{4} = \frac{n}{4}. \nonumber
\end{equation}
Suppose that a component $A$ of $G - U_x^+$ is a triangle. Then, $V(A) = N[v] \backslash \{x\}$. Let $Y = U_x^+ \cup V(A)$. Since $G-Y$ contains no triangle, $\iota_{\rm c}(G-Y) \leq \frac{n-|Y|}{4}$ by Lemma~\ref{lemmacomp} and the induction hypothesis. Let $D_{G-Y}$ be a cycle isolating set of $G-Y$ of size $\iota_{\rm c}(G-Y)$. Since $v \in N(x) \cap V(A)$ and $U_x^+ \subset N[x]$, $\{x\} \cup D_{G-Y}$ is a cycle isolating set of $G$. Thus, $\iota_{\rm c}(G) \leq 1 + \frac{n-|Y|}{4}$. Since $|Y| \geq 7$, $\iota_{\rm c}(G) < \frac{n}{4}$. 

Now suppose $|U_x| \leq 2$. Then, $1 \leq |\mathcal{H}_x'| \leq 2$.\medskip 
\\
\textit{Case 1: $|\mathcal{H}_x'| = 1$.} Let $T$ be the member of $\mathcal{H}_x'$. We have $xy \in E(G)$ for some $y \in V(T)$. Let $Y = \{x\} \cup V(T)$. Then, $Y \subseteq N[y]$. Also, $G-Y$ has a component $T'$ with $N[v] \backslash \{x\} \subseteq V(T')$. Since $T$ is the only member of $\mathcal{H}'$ that is linked to $x$, the components of $G-Y$ are $T'$ and the members of $\mathcal{H}_x^*$. Recall that no member of $\mathcal{H}_x^*$ is a triangle.

If $T'$ is not a triangle, then, by Lemma~\ref{lemma} (with $X = \{y\}$), Lemma~\ref{lemmacomp}, and the induction hypothesis, we have $\iota_{\rm c}(G) \leq 1 + \iota_{\rm c}(G - Y) \leq 1 + \frac{n-|Y|}{4} = \frac{n}{4}$. 

Suppose that $T'$ is a triangle. Let $W = V(T) \cup V(T')$. We have $x \notin W$, $y \in N(x) \cap V(T)$, and $v \in N(x) \cap V(T')$. Also, the components of $G-\{x\}$ are the components of $G[W]$ and the members of $\mathcal{H}_x^*$. By the induction hypothesis, each member $H$ of $\mathcal{H}_x^*$ has a cycle isolating set $D_H$ with $|D_H| \leq \frac{|V(H)|}{4}$.

Suppose that $G[W]-(N(x) \cap W)$ contains no cycle. Then, $\{x\} \cup \bigcup_{H \in \mathcal{H}_x^*} D_H$ is a cycle isolating set of $G$. Thus, $\iota_{\rm c}(G) \leq 1 + \sum_{H \in \mathcal{H}_x^*} \frac{|V(H)|}{4} = 1 + \frac{n-|\{x\} \cup W|}{4} = 1 + \frac{n-7}{4} < \frac{n}{4}$.

Now suppose that $G[W]-(N(x) \cap W)$ contains a cycle $A$. Since $v, y \in N(x) \cap W$ and $|W| = 6$, either $A \simeq C_3$ or $A \simeq C_4$. 

Suppose $A \simeq C_4$. Then, $N(x) \cap W = \{v, y\}$, $V(A) = (V(T) \backslash \{y\}) \cup (V(T') \backslash \{v\})$, and hence 
$uw \in E(A) \subseteq E(G)$ for some $u \in V(T) \backslash \{y\}$ and some $w \in V(T') \backslash \{v\}$. Let $Z = \{w\} \cup V(T)$ and let $x'$ be the member of $V(T') \backslash \{v,w\}$. We have $V(G-Z) = \{v, x, x'\} \cup \bigcup_{H \in \mathcal{H}_x^*} V(H)$ and $x, x' \in N(v)$. Thus, since the members of $\mathcal{H}_x^*$ are linked to $x$, $G-Z$ is connected. Since $N(x) \cap W = \{v, y\}$, we have $xx' \notin E(G-Z)$, so $G-Z$ is not a triangle. By the induction hypothesis, $\iota_{\rm c}(G-Z) \leq \frac{n-|Z|}{4} = \frac{n-4}{4}$. Since $Z \subseteq N[u]$, Lemma~\ref{lemma} (with $X = \{u\}$) gives us $\iota_{\rm c}(G) \leq 1 + \iota_{\rm c}(G-Z) \leq \frac{n}{4}$.

Now suppose $A \simeq C_3$. Since $V(A) \subseteq W \backslash N(x) \subseteq (V(T) \cup V(T')) \backslash \{v, y\}$, $V(A)$ contains either the two vertices in $V(T) \backslash \{y\}$ and one of the two vertices in $V(T') \backslash \{v\}$ or the two vertices in $V(T') \backslash \{v\}$ and one of the two vertices in $V(T) \backslash \{y\}$. Suppose $|V(A) \cap (V(T') \backslash \{v\})| = 1$. Then, $V(T) \backslash \{y\} \subseteq V(A)$. Let $x'$ be the member of $V(A) \cap (V(T') \backslash \{v\})$. Let $Z = (V(T) \backslash \{y\}) \cup V(T')$. Since $V(A) \cup V(T') \subseteq N[x']$, $Z \subseteq N[x']$. We have $V(G-Z) = \{x, y\} \cup \bigcup_{H \in \mathcal{H}_x^*} V(H)$. Since $xy \in E(G)$ and the members of $\mathcal{H}_x^*$ are linked to $x$, $G-Z$ is connected. Since $N(y) \cap \left( \bigcup_{H \in \mathcal{H}_x^*} V(H) \right) = \emptyset$, $G-Z$ is not a triangle. By the induction hypothesis, $\iota_{\rm c}(G-Z) \leq \frac{n-|Z|}{4} = \frac{n-5}{4}$. Since $Z \subseteq N[x']$, Lemma~\ref{lemma} (with $X = \{x'\}$) gives us $\iota_{\rm c}(G) \leq 1 + \iota_{\rm c}(G-Z) < \frac{n}{4}$. Similarly, $\iota_{\rm c}(G) < \frac{n}{4}$ if $|V(A) \cap (V(T) \backslash \{y\})| = 1$.\medskip
\\
\textit{Case 2: $|\mathcal{H}_x'| = 2$.} Let $T_1$ and $T_2$ be the two members of $\mathcal{H}_x'$. 

Suppose that $T_2$ is linked to a member of $L \backslash \{x\}$. Since $T_1$ is linked to $x$, $xy \in E(G)$ for some $y \in V(T_1)$. Let $Y = \{x\} \cup V(T_1)$. Then, no component of $G-Y$ is a triangle ($G-Y$ has a component $A$ with $(N[v] \backslash \{x\}) \cup V(T_2) \cup \bigcup_{H \in \mathcal{H}' \backslash \{T_1, T_2\}} V(H) \subseteq V(A)$, and the other components of $G-Y$ are the members of $\mathcal{H}_x^*$). Also, $Y \subseteq N[y]$. By Lemma~\ref{lemma} (with $X = \{y\}$), Lemma~\ref{lemmacomp}, and the induction hypothesis, we have 
\begin{equation} \iota_{\rm c}(G) \leq 1 + \iota_{\rm c}(G - Y) = \frac{|Y|}{4} + \sum_{H \in {\rm C}(G - Y)} \iota_{\rm c}(H) \leq \frac{|Y|}{4} + \sum_{H \in {\rm C}(G - Y)} \frac{|V(H)|}{4} = \frac{n}{4}. \nonumber
\end{equation}
Similarly, $\iota_{\rm c}(G) \leq \frac{n}{4}$ if $T_1$ is linked to a member of $L \backslash \{x\}$. 

Now suppose that, for each $i \in \{1, 2\}$ and each $x' \in L \backslash \{x\}$, $T_i$ is not linked to $x'$. Then, $T_1$ and $T_2$ are components of $G-\{x\}$. Let $Y = \{x\} \cup V(T_1) \cup V(T_2)$. 

Suppose that no component of $G - \{x\}$ other than $T_1$ and $T_2$ is a triangle. Then, $\iota_{\rm c}(G-Y) \leq \frac{n-|Y|}{4}$ by Lemma~\ref{lemmacomp} and the induction hypothesis. Let $D_{G-Y}$ be a cycle isolating set of $G-Y$ of size $\iota_{\rm c}(G-Y)$. Since $N(x) \cap V(T_1) \neq \emptyset \neq N(x) \cap V(T_2)$, $\{x\} \cup D_{G-Y}$ is a cycle isolating set of $G$, so $\iota_{\rm c}(G) \leq 1 + \frac{n-|Y|}{4} = 1 + \frac{n-7}{4} < \frac{n}{4}$. 

Now suppose that $G - \{x\}$ has a component $T_3$ such that $T_3 \notin \{T_1, T_2\}$ and $T_3$ is a triangle. Since $T_1$ and $T_2$ are the only members of $\mathcal{H}'$ that are linked to $x$, it follows that $V(T_3) = N[v] \backslash \{x\}$. Let $Y' = Y \cup V(T_3)$. Since $G-Y'$ contains no triangle, $\iota_{\rm c}(G-Y') \leq \frac{n-|Y'|}{4}$ by Lemma~\ref{lemmacomp} and the induction hypothesis. Let $D_{G-Y'}$ be a cycle isolating set of $G-Y'$ of size $\iota_{\rm c}(G-Y')$. Since $v \in N(x) \cap V(T_3)$ and $N(x) \cap V(T_1) \neq \emptyset \neq N(x) \cap V(T_2)$, $\{x\} \cup D_{G-Y'}$ is a cycle isolating set of $G$. Thus, $\iota_{\rm c}(G) \leq 1 + \frac{n-|Y'|}{4} = 1 + \frac{n-10}{4} < \frac{n}{4}$.~\hfill{$\Box$}

\end{document}